\documentclass {article}
\usepackage[left=1.75in, right=1.75in, top=1in, bottom=1in]{geometry}

\usepackage{titling}
\usepackage{titlesec}
\usepackage{lipsum}
\usepackage{cite}
\usepackage{float}
\usepackage[toc,page]{appendix}
\usepackage{graphicx}
\usepackage{listings}
\usepackage{color}
\usepackage[fleqn]{amsmath}
\usepackage{hyperref}

\definecolor{dkgreen}{rgb}{0,0.6,0}
\definecolor{gray}{rgb}{0.5,0.5,0.5}
\definecolor{mauve}{rgb}{0.58,0,0.82}
\hypersetup{hidelinks}

\pretitle{\begin{center}\LARGE\bfseries}
\posttitle{\par\end{center}\vskip 0.5em}
\preauthor{\begin{center}\Large\ttfamily}
\postauthor{\end{center}}
\predate{\par\large\centering}
\postdate{\par}
\titleformat{\section}[block]{\Large\bfseries\filcenter}{}{1em}{}
\titleformat{\subsection}[hang]{\bfseries}{}{1em}{}
\title{A Survey of Mathematical Models on Somitogenesis}
\author{Hanyu Song\\*Brandeis University}
\date{September 09, 2019}

\begin{document}
\maketitle
\pagenumbering{arabic}
\thispagestyle{plain}

\begin{abstract}
This paper presents a comprehensive survey of various established mathematical models pertaining to Somitogenesis, a biological process. The study begins by revisiting and replicating the findings from prominent research papers in this domain, subsequently offering a critical evaluation of the strengths and weaknesses inherent in each approach. By synthesizing this knowledge, the paper aims to contribute to a deeper understanding of Somitogenesis, and pave the way for further advancements in the development of enhanced mathematical models for this intricate biological process. The concluding section offers valuable insights and directions for prospective research in this field.
\end{abstract}

\tableofcontents
\newpage
 
\section{Introduction to Somitogenesis}
\label{sec: intro}
Somites are blocks of cells that lie along the anterior-posterior (AP) vertebrate embryonic axis of the developing embryo.
\par
\begin{figure}[H]
  \centering
  \includegraphics[width=\linewidth]{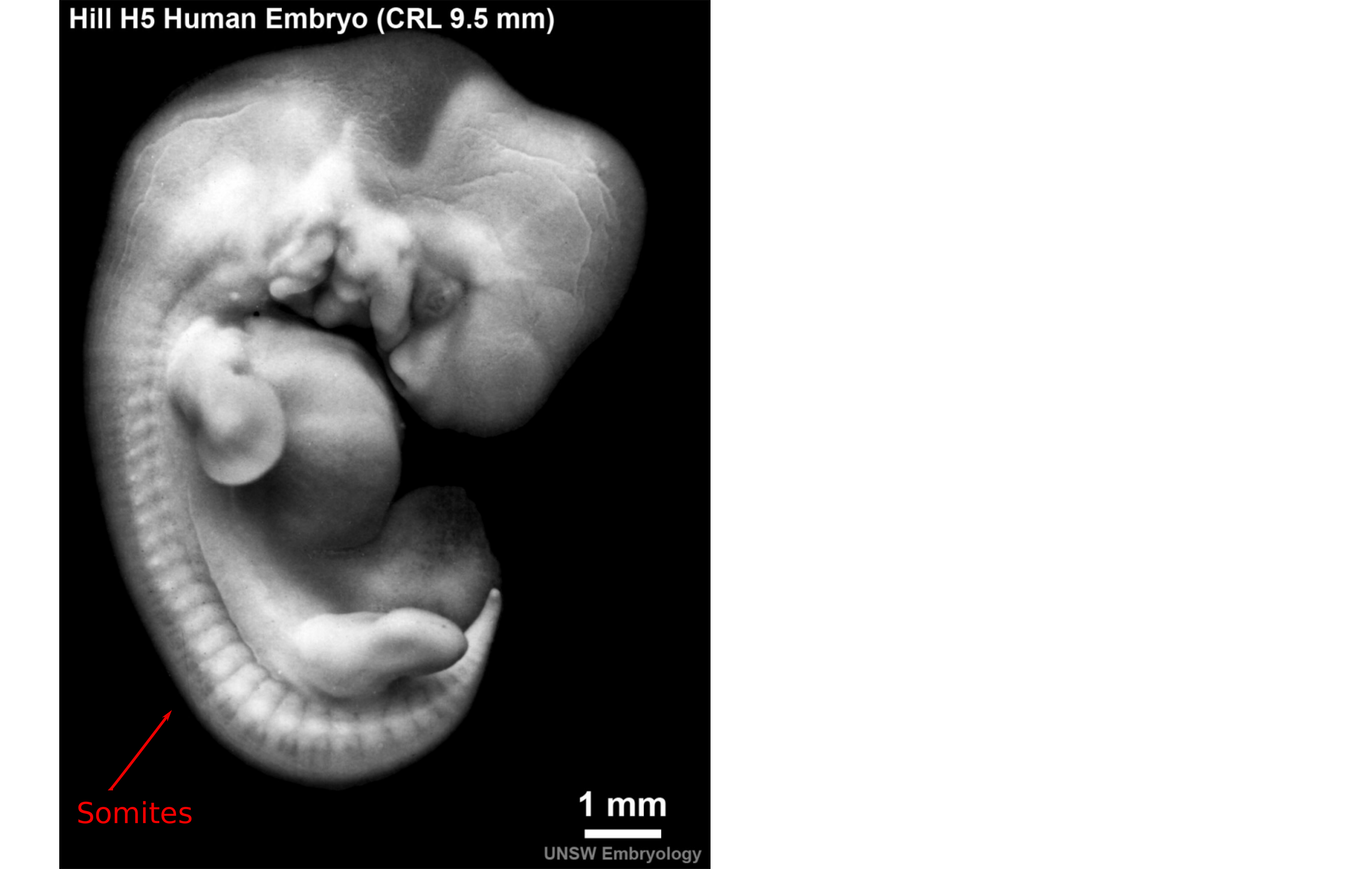}
  \caption{Embryonic somite and the AP axis. The left picture is a human embryo\cite{stage16}, where the somites are already in shape when the embryo is still very immature. The right picture is an anterior-posterior axis where somites segment from the PSM cell that lies on both sides of the AP axis\cite{Baker2006}. Front the posterior end to the anterior, the cells transform from undetermined, to determined, then finally to somites.}
  \label{fig: intro}
\end{figure}

Somitogenesis is the process by which somites form by segmenting the axis into similar morphological units such as vertebrates etc... Somitogenesis serves as the key biological process in the embryo since it's responsible for segmenting the vertebrate axis and generating the prepattern that guides the formation of the tendons, ribs, muscle, and other associated features of the body trunk. Figure \ref{fig: intro} illustrates the form of somites in an embryo and how segmentation works in the AP axis. 

\par

\begin{figure}[H]
  \centering
  \includegraphics[width=\linewidth]{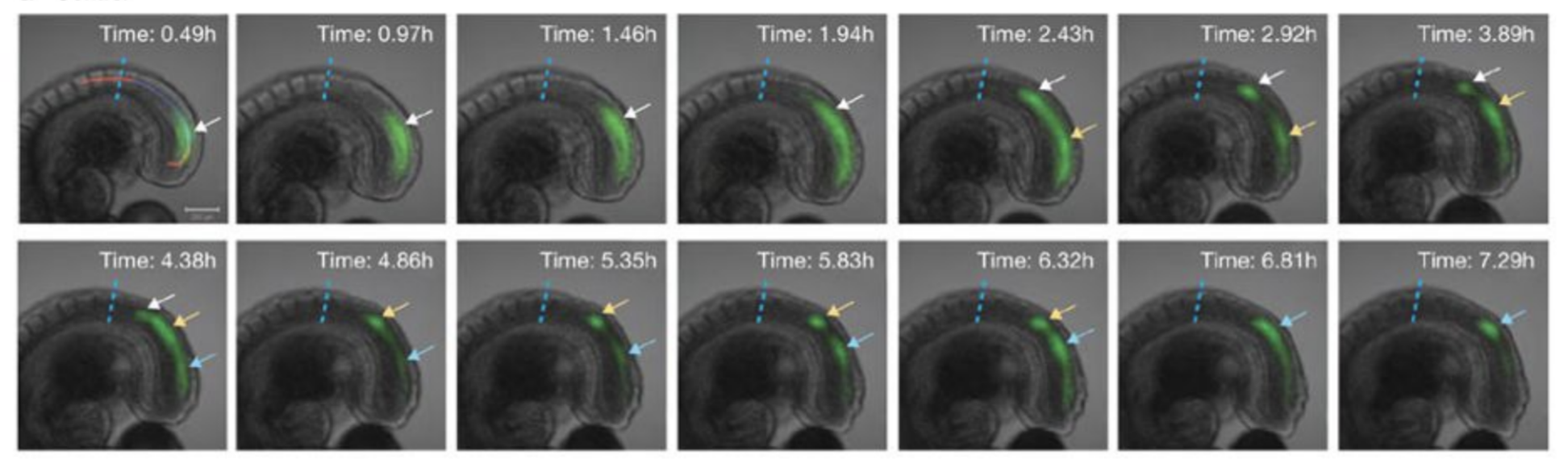}
  \caption{Gene expression in a mouse embryo, where the gene is marked green\cite{Aul}. The gene has a wave-like propagation from the posterior end of the PSM to the anterior side.}
  \label{fig: fgf8wave}
\end{figure}

Although many details about somitogenesis are still debated, there are some scientific facts that serve as the fundamentals for further research: Somites segment from the presomitic mesoderm (PSM): thick bands of tissue that lie on either side of the AP axis. The segmentation begins with the establishment of a prepattern of gene expression, and it is characterized by periodic activations in regions where future somites will segment. Early scanning microscope images show that the posterior PSM displays a series of cells similar in size and structure, known as \textit{somitomeres}, which seem to be the precursors of the somites\cite{Gossler}. The existence of this prepattern was confirmed by microsurgical experiments in which isolated parts of the PSM formed somites in strict isolation\cite{Chernoff}. Figure \ref{fig: fgf8wave} demonstrates the wave-like gene expression in a mouse embryo. 

\par

Another fact about the PSM is that the PSM is not a homogeneous tissue\cite{Merier}. This is supported by microsurgical experiments conducted by Dubrulle and co-workers: AP inversions of somite-length regions of the posterior PSM resulted in normal segmentation whilst inversions of the anterior PSM resulted in somites with reversed polarity\cite{Dubrulle2001}, which suggested that the anterior-most part of the PSM is determined with regard to its segmentation program, whilst the posterior-most part of PSM is susceptible in this respect. This proves the PSM's heterogeneity, which is a key feature of the models for somitogenesis. 
\par
\begin{figure}[H]
  \centering
  \includegraphics[width=\linewidth]{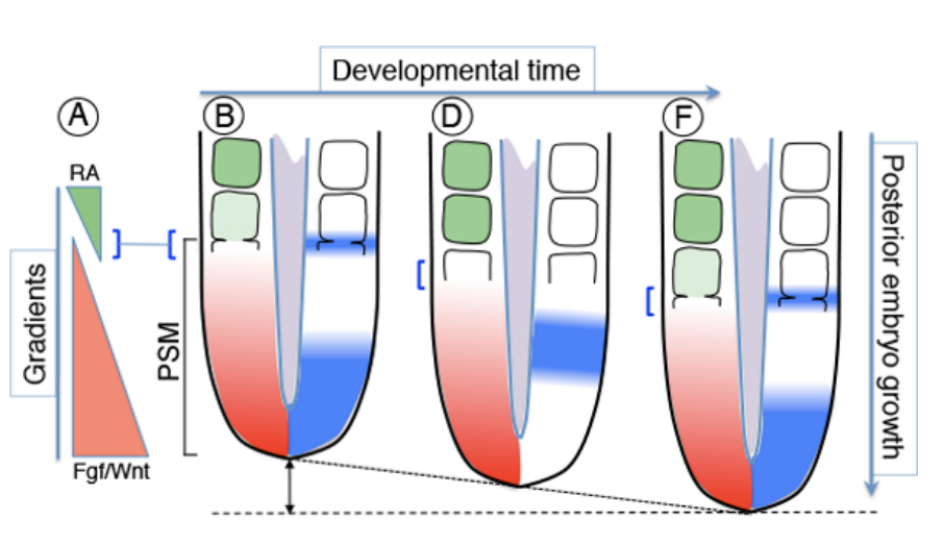}
  \caption{PSM cell and gene gradients, FGF included\cite{Mallo}. The PSM elongates posteriorly as the somites are formed, whilst the gradient of FGF8 always has a peak at the very posterior end of the PSM cell. As it decreases to a certain level, the prepattern arises (the blocks that are not fully green), then somites are formed. What's contrasting to the FGF gradient is the RA gradient, which is another type of gene that peaks at the anterior end of the PSM, but it's not that relevant to the overall process compared to FGF8.}
  \label{fig: PSM}
\end{figure}

The different regions of the PSM were found to correspond to regions of varying FGF signaling. \textit{fgf8}, which is a gradient of FGF8 (Fibroblast Growth Factor 8). FGF8 is a gene with dynamic expression in the PSM, peaking at the posterior end of the embryo, whilst decreasing in the direction of the anterior end\cite{Dubrulle2004a}. See Figure \ref{fig: PSM}. The function of \textit{fgf8} is to down-regulate the cells, meaning, higher concentration of \textit{fgf8} will prevent the segmentation of PSM, whilst its decrease will make segmentation possible, and when \textit{fgf8} decreases past a certain threshold, the cells are then able to segment into somites. We call that threshold "the determination front"\cite{Dubrulle2001}. The uneven distribution of \textit{fgf8} implies that the positional information of the PSM cell is crucial. However, the role of positional information is a controversial issue in mathematical biology, and it's typically not possible to build robust biological structures without additional mechanisms, such as diffusion\cite{Touboul}.
\subsection{Unsolved questions}
We know the important information that the down-regulation of \textit{fgf8} heavily affects the somitogenesis process, and we seem to understand the logic behind the process of somitogenesis, but it is difficult to draw any conclusion about which specific type of model is capable of accurately recapitulating this process. There are still many questions that must be determined: whether the PSM cells are oscillatory or excitable with respect to \textit{fgf8} levels? Are the cells globally controlled by the gene or do they have local interactions between themselves as well? Does the global \textit{fgf8} down-regulation even matter? What will happen if the \textit{fgf8} is kept constant, can a reaction-diffusion model, that emphasizes local interactions between cells, explain this process accurately? In the later discussions of this paper, we will look into several kinds of mathematical models, with each of them having distinct answers to the above questions. Admittedly, none of them are deemed to be "perfect", with each of them having its own drawbacks. Although there are plenty of models out there, understanding those models' mechanisms is still important as it could accelerate the process of creating a better one in the future.
\par
Although there are plenty of models out there, to date nobody has provided any comparison of them, and most papers on this topic don't even reference each other as they are in different fields: Mathematical Biology, Bio Development, Physics, etc. Therefore, this paper's goal is not to compare and find the perfect model, but to see each of their distinctive advantages and try to synthesize them if possible, while avoiding their drawbacks when we attempt to create new models in the future. 
\section{Clock and Wavefront Model}
\label{sec: clock and wavefront}
\subsection{Summary}
One of the most famous and widely studied models is the clock and wavefront (C \& W) model. As its name implies, the model proposes the existence of a segmentation clock and a wavefront of FGF8 along the AP axis of vertebrate embryos. This idea was first proposed in 1975 by Cooke and Zeeman, with the gist that there is a longitudinal, global positional information, which is the above-mentioned FGF8 gradient distribution, that interacts with a smooth cellular oscillator, which is the so-called clock, to govern the time for the PSM cells to segment and develop into somites. This idea was then revised by Pourquie and co-workers, where they went into more specifics and proposed that the clock sets the times at which new somite boundaries form whilst the position of the determination front sets where they form\cite{Dubrulle2001}. For a cell at a particular point, they assume that competence to segment is only achieved once FGF8 signaling has decreased below a certain threshold, the position of which is known as the determination front. 

\begin{figure}[H]
  \centering
  \includegraphics[scale = 0.75]{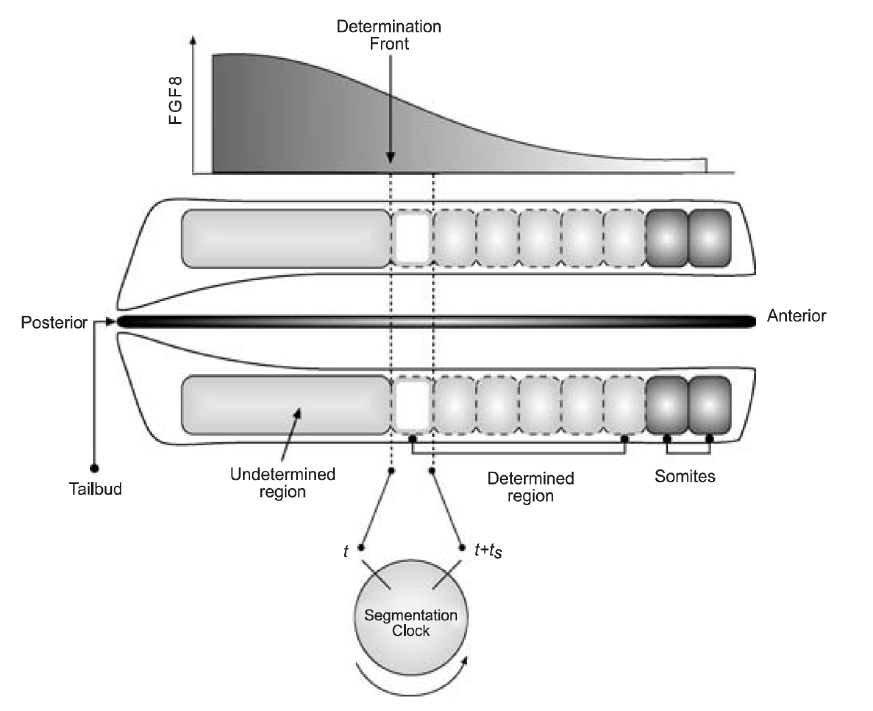}
  \caption{Representation of the vertebrate body plan during somite formation\cite{Baker2006b}. The top part of the diagram shows the FGF8 wavefront, with a peak in the posterior and a decrease in the direction of the anterior. When the FGF8 decreases to a certain level, it reaches the determination front. The middle section of the diagram shows the AP axis of the embryo with the somites (dark grey blocks), determined region (light grey blocks), and the undetermined region (light grey band) clearly marked. The bottom is a visualization of a segmentation clock, which shows the time needed for cells to gain segmentation.}
  \label{fig: clockandwavefront}
\end{figure}

\par
Therefore, the whole somitogenesis process, according to this model, is divided and analyzed into different parts. Before reaching a determination front, a cell will gain the ability to segment by being able to produce a "somitic factor", which could be several possibilities of genes. One clock oscillation after reaching the determination front, cells become able to produce the "signaling molecule". After a cell is able to produce somitic factor and respond to the signaling molecule, it is specified as somitic and becomes refractory to FGF8 signaling\cite{Baker2008} 

\subsection{Mathematical Equations}
C \& W mathematical equations were first proposed by Collier \textit{et al.} (2000) and were developed by Mclnerny \textit{et al.} in 2004, then by Baker and colleagues in 2006. One of the most important features of this model is that local mechanisms, controlled by time points and positional information, will trigger segmentation, which fits the C \& W assumption perfectly. After segmentation, cells adhere to each other, creating distinctive somites. When creating this model, Collier made some further assumptions\cite{McInerney2004}: (1) The AP-axis can be seen as fixed with respect to the cells. The PSM's length is constant and the segmentation pattern progresses with a constant speed. (In reality, the posterior end is actually elongating.) (2) The signals that are emitted by specified cells when they reach certain points are like pulses. The signaling molecule disperses fast and diffuses rapidly. This is the key assumption since rapid diffusion can ensure that only cells that are in certain positions will respond to the signal, if not, all cells will segment at the same time. (3) Somites are formed continually, and the beginning or end of this process is not considered, which means it's assumed that signals emitted from cells exist all the time. 
\par
\begin{figure}
  \centering
  \includegraphics[scale = 0.42]{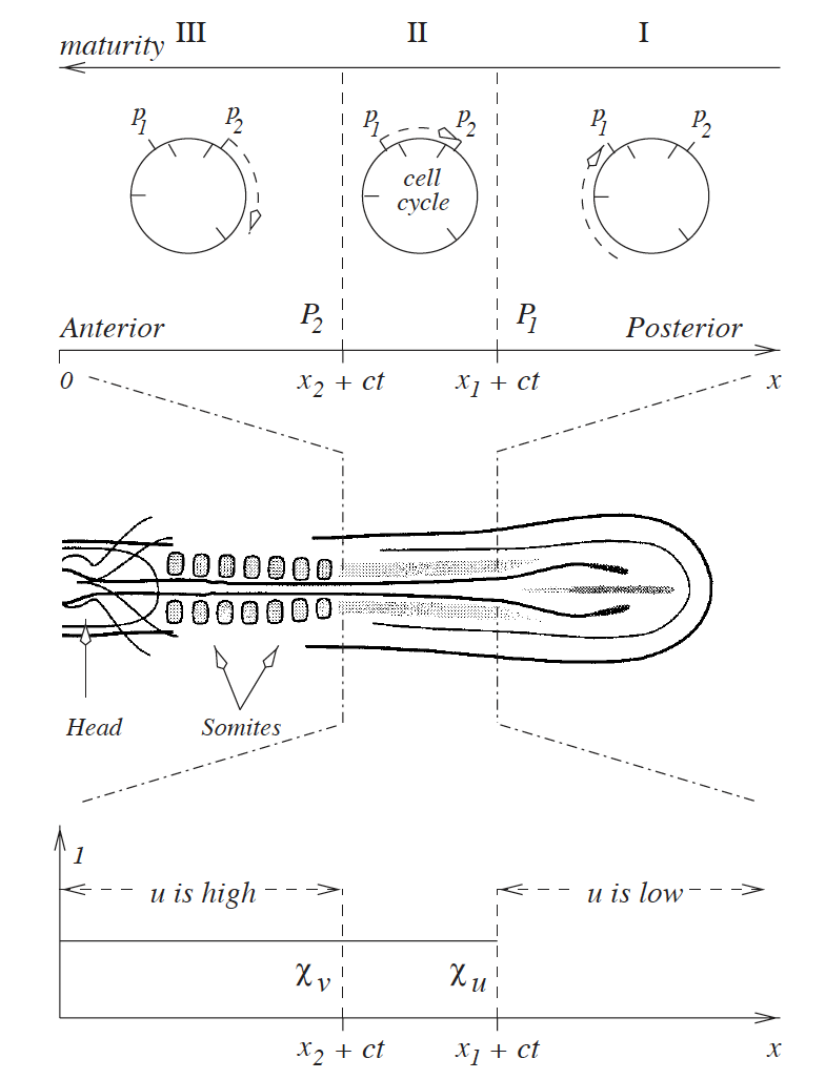}
  \caption{Representation of C \& W model illustrating the two time points P1 and P2 and the three key stages within the model.\cite{McInerney2004}. Cells at the posterior end of the PSM (Region \uppercase\expandafter{\romannumeral1}) are less mature than those in other regions since the somitogenesis process starts from the anterior part. As cells become more mature in Region  \uppercase\expandafter{\romannumeral2}, they become capable of responding to the signaling molecule v, emitted by cells at point P2. In Region  \uppercase\expandafter{\romannumeral3} they begin to form somites and are no longer able to emit any signals.}
  \label{fig: cwds1}
\end{figure}

 This model can be well explained by Figure \ref{fig: cwds1}. In this diagram, x denotes the distance while t denotes time. There are two key components in this model: u(x,t) and v(x,t), in which u(x,t) represents the degree of concentration of somitic factor a cell is exposed to at a given x and t, while v(x,t) represents the diffusive signaling molecule. A cell that has a high concentration of u is specified as somitic, while those with a low concentration of u are non-somitic. 

\par The mathematical equations proposed by Collier and colleagues are also based on these two components \cite{McInerney2004}: 
\begin{gather}
\partial_t u(x,t) = \frac{(u+\mu v)^2}{\gamma+\kappa u^2}\chi_u(x,t)-\frac{u}{k} \\
\partial_t v(x,t) = \frac{\chi_v(x,t)}{\epsilon+u}-v-D\frac{\partial^2 v}{\partial x^2}
\end{gather}
$\chi_u$ and $\chi_v$ are controlled by two Heaviside step functions: 
\begin{gather}
\chi_u = H(ct-x+x_1) \\
\chi_v = H(ct-x+x_2)
\end{gather}
note that the Heaviside step functions' rule of calculation is: 
\begin{gather}
H(x)=\left\{
\begin{aligned}
1 &  & x\geq0 \\
0 &  & x<0
\end{aligned}
\right.
\end{gather}
\par As mentioned above, u and v represent the concentration of "somitic growth factor" and "signaling molecule" respectively, while other variables in these equations are all positive constants. This model uses a zero flux boundary condition. It indicates that this boundary condition prevents anything from leaving this system, which may be an application of the third assumption made by Collier mentioned above. 
\par Heaviside functions,$\chi_u$ and $\chi_v$, play an important role in this model. They can be seen as switches: the elements inside brackets, t, and x, which are the time and location information, together determine the on and off of the dynamics in u and v respectively. In Figure \ref{fig: cwds1}, the Heaviside functions are shown along with the regions where the somitic growth factor u is, respectively, high ($x < x_2+ct$) and low ($x > x_1+ct$). Somitic growth factor and signaling molecules boost the somitogenesis process collectively and they affect each other, as we can see that $\partial_t u$ is affected by v and $\partial_t v$ is affected by u. Specifically, u inhibits v while v activates u. 
\par The model was further expanded by Baker and colleagues in 2006. His team made some revisions to the two equations above and they added a third equation into the system: $\frac{\partial w}{\partial t}$ which represents the changing gradient of FGF8, which down-regulates the somitogenesis process\cite{Baker2006b}: 
\begin{gather}
\frac{\partial u}{\partial t} = \frac{(u+\mu v)^2}{\gamma+u^2}\chi_u-u \\
\frac{\partial v}{\partial t} = k(\frac{\chi_v}{\epsilon+u}-v)+D_v\frac{\partial^2 v}{\partial x^2} \\
\frac{\partial w}{\partial t} =\chi_w-\eta w+D_w\frac{\partial^2 w}{\partial x^2}
\end{gather}
and $\chi_w = H(x-x_n-c_n t)$ where $x_n$ and $c_n$ are constants. Based on the previous system, this system is reproducing these important behaviors: (1) the increase of somitic factor u is activated by signaling molecules and is self-regulating. (2) The somitic factor is an inhibiting signaling molecule. In other words, signaling molecule is produced rapidly in areas where somitic factor concentration is low. (3) FGF8 is produced in the tail and regresses along the x-axis. 

\subsection{Analysis}
The C \& W mathematical model proved to be effective in producing a qualitatively reasonable match to reality. We recapitulated and reproduced some of the results of the above mathematical equations shown in Mclnerney and Baker's papers, and they do support the gist of C \& W theory. We first analyzed the qualitative behaviors of this model, and the result can be explained with Figure \ref{fig: qualitative}, which is derived from Figure \ref{fig: cwds1}:

\begin{figure}[H]
  \centering
  \includegraphics[width = \linewidth]{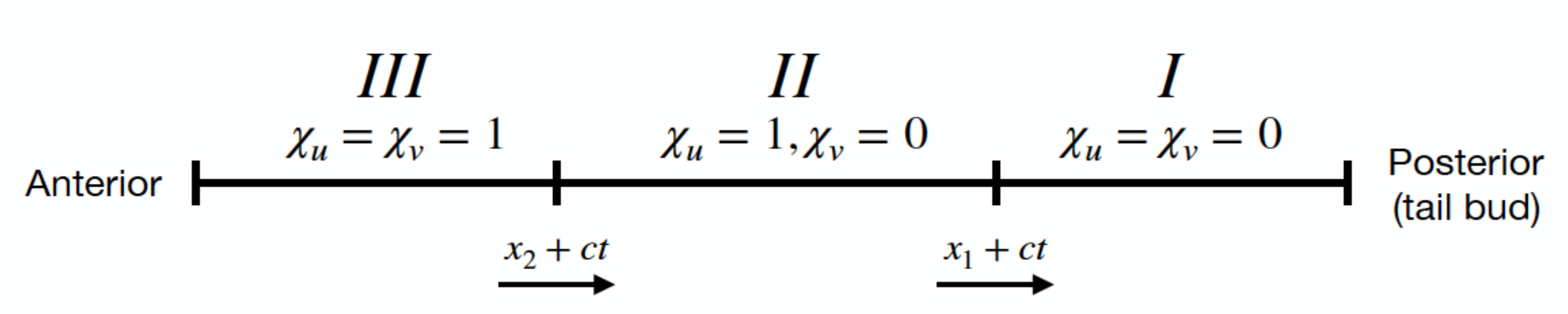}
  \caption{The three stages of somitogenesis from the posterior to the anterior end of the PSM.}
  \label{fig: qualitative}
\end{figure}

In region \uppercase\expandafter{\romannumeral1}, since the switches $\chi_u = \chi_v = 0, (u,v) \rightarrow (0,0)$, while in Region \uppercase\expandafter{\romannumeral2} and \uppercase\expandafter{\romannumeral3}, as $ \chi_u and \chi_v $ change with respect to t, the qualitative behavior of these two regions will be different. Below are the phase planes of u and v in Region \uppercase\expandafter{\romannumeral2} and \uppercase\expandafter{\romannumeral3} respectively, which is a reproduction of Fig 2 and Fig 3 in Mclnerney's paper\cite{McInerney2004} using XPP: 
\begin{figure} [H]
  \centering
  \includegraphics[width = \linewidth]{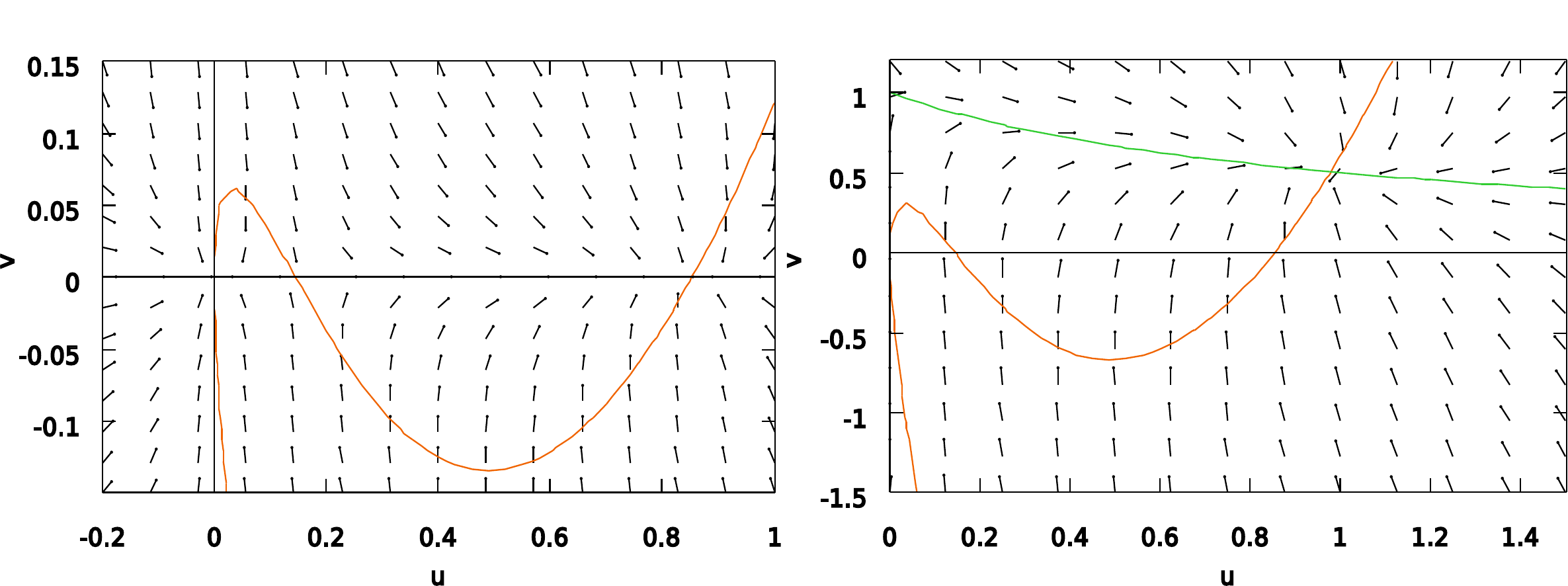}
  \caption{The phase planes of u and v in Region \uppercase\expandafter{\romannumeral2} (left) and Region \uppercase\expandafter{\romannumeral3} (right).}
  \label{fig: xpp}
\end{figure}
In Region \uppercase\expandafter{\romannumeral2}, there are three steady states: two stable equilibriums when u is close to 0 and 1, and a saddle in the middle. The region \uppercase\expandafter{\romannumeral2} carries the pulse of the signaling molecule. In the phase plane, we can see that after the cells pass the determination front and before they finish one clock cycle, the cells gain the ability to respond to signals.- the somitic factor concentration u is always above 0, while they can't generate signaling molecules themselves, as v remains to be 0. In Region \uppercase\expandafter{\romannumeral3}, after cells undergo one cycle of segmentation clock, $\chi_u = \chi_v=1$. There is only one stable equilibrium in the phase plane, meaning whatever u and v's initial values are, they will all arrive at that specific point and their values will remain stable, then cells are identified as somitic. 
\par
\begin{figure} [H]
  \centering
  \includegraphics[width = \linewidth]{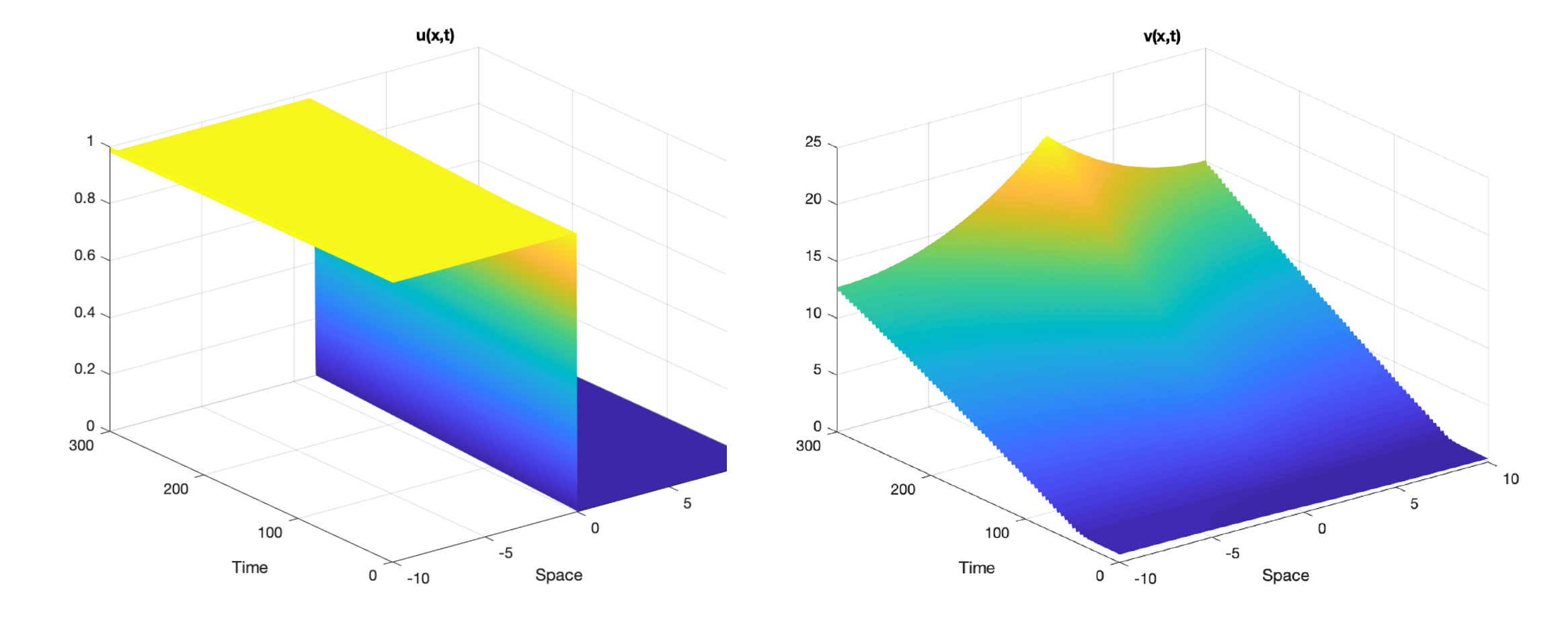}
  \caption{Numerical solution given by equation (1) and (2) for $0 \leq t \leq 300$. Parameter values: $\mu = 10^{-1}, \gamma = 0.2, \kappa = 10, c = 5*10^{-3}, \epsilon = 10^{-3}, D = 100$. This set of parameter values violated one of the conditions mentioned in the paper, so this solution is not a perfect one since for $\mu = 10^{-5}$ and $\gamma = 0.2$, a high level of v fails to activate u production\cite{McInerney2004}.}
  \label{fig: solution1}
\end{figure}

 We then reproduced the numerical solutions for equations (1) and (2) in Mclnerney's Figure 11(a) and (b)\cite{McInerney2004}, shown in figure \ref{fig: solution1}.

We also reproduced the numerical solution of the C \& W model in one spatial dimension given by equations (5) (6) (7 )in Baker's Fig 3\cite{Baker2008}, using the code provided in the appendix. Figure \ref{fig: solution2} contains the numerical solutions for $u(x,t), v(x,t)$ and $w(x,t)$ respectively: 

\begin{figure} [H]
  \centering
  \includegraphics[width = \linewidth]{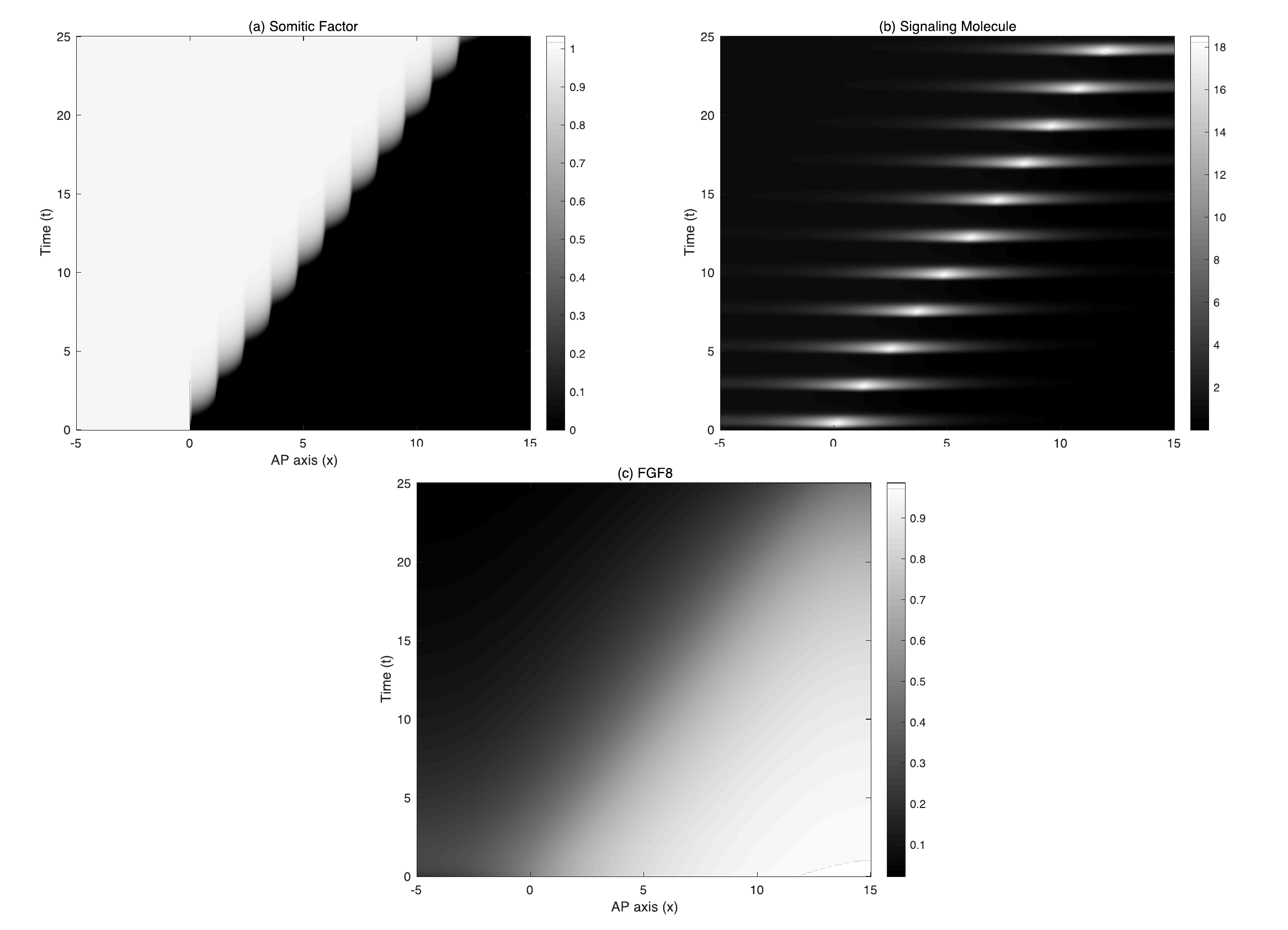}
  \caption{Numerical solution given by equation (5) (6) and (7), showing the spatiotemporal dynamics of the somitic factor, (a). the signaling molecule, (b), and \textit{fgf8} (c). The regression of the FGF8 wavefront is accompanied by a series of pulses in the signaling molecule and coherent rises in the level of the somitic factor. Parameter values:  $\mu = 10^{-4}, \gamma = 10^{-3}, \kappa = 10, \epsilon = 10^{-3}, \eta = 1.0, D_v = 50, D_w = 20, x_n = 0, c_n = 0.5, D = 100$.}
  \label{fig: solution2}
\end{figure}
\par However, the fact that the verification of the results of the above models, in some ways, shows the models' validity, can not prove the models to be flawless. There are some issues that need to be considered before constructing a better model. The equation set (1) and (2) is not robust because somites depend sensitively on plenty of factors, such as mesh, the speed c, and initial conditions. Any slight interference in those factors will prevent them from obtaining successful results. For equations (5), (6), and (7), although the results in Fig \ref{fig: solution2} show a clear and consistent pattern of pulses of somitic factor and signaling molecule, they rely heavily on a very smooth gradient w. Admittedly, in normal cases, the idea that u and v rely on a smooth gene gradient is not, in itself, problematic. However, in this model, it's simply assumed that a \textit{generic} FGF8 molecule makes up the gradient controlling the position of the determination front\cite{Baker2008}, which means although the gene's name is FGF8, it in fact represents the aggregate influence of multiple genes that may affect the somitogenesis process. In other words, the gradient is modeled at a very phenomenological level. Requiring such a gradient to be perfectly smooth becomes a drawback of this system. A stochastic FGF8 gradient or some random \textit{\textit{fgf8}} pulse will easily mess up the result. This problem is mentioned and demonstrated in Fig 3 in Baker's paper\cite{Baker2008}. Also, in this paper, the position of the determination front is prescribed, yet in reality, it will be subject to plenty of factors such as the gradient slope, etc... 

\section{Oscillatory-based Model}
\label{sec: PORD}
\subsection{Summary of the PORD Model}
\par In the C \& W model mentioned above, the key is that long-range molecular gradients control the movement of the front and therefore the placement of the stripes in the embryo. In this section, we are introducing a fundamentally different system: the progressive oscillatory reaction-diffusion model, or PORD model, which does not rely on a global gradient control, but is driven by short-range interactions. 
\par In the first section of this paper, we have introduced several "facts" and "unsolved questions" in the somitogenesis process. Although the oscillatory model's mechanism is very different from C \& W, they do share a lot of similarities - it's just that their interpretations and understandings to those facts are different. The PORD model admits the existence of the posterior movement of the determination front, yet it explains that it's not controlled by global positional information but by interactions between cells. In Cotterell's paper, it's argued that the PORD model could also explain some other important features of somitogenesis, such as size regulation, which previous reaction-diffusion models fail to explain. However, we did find that controlling the FGF8 gradient such as adding random pulse in the C \& W model, will result in larger somites, which strengthens the argument that the amount of somitic factor controls the size of somites. 

\begin{figure} [H]
  \centering
  \includegraphics[width = \linewidth]{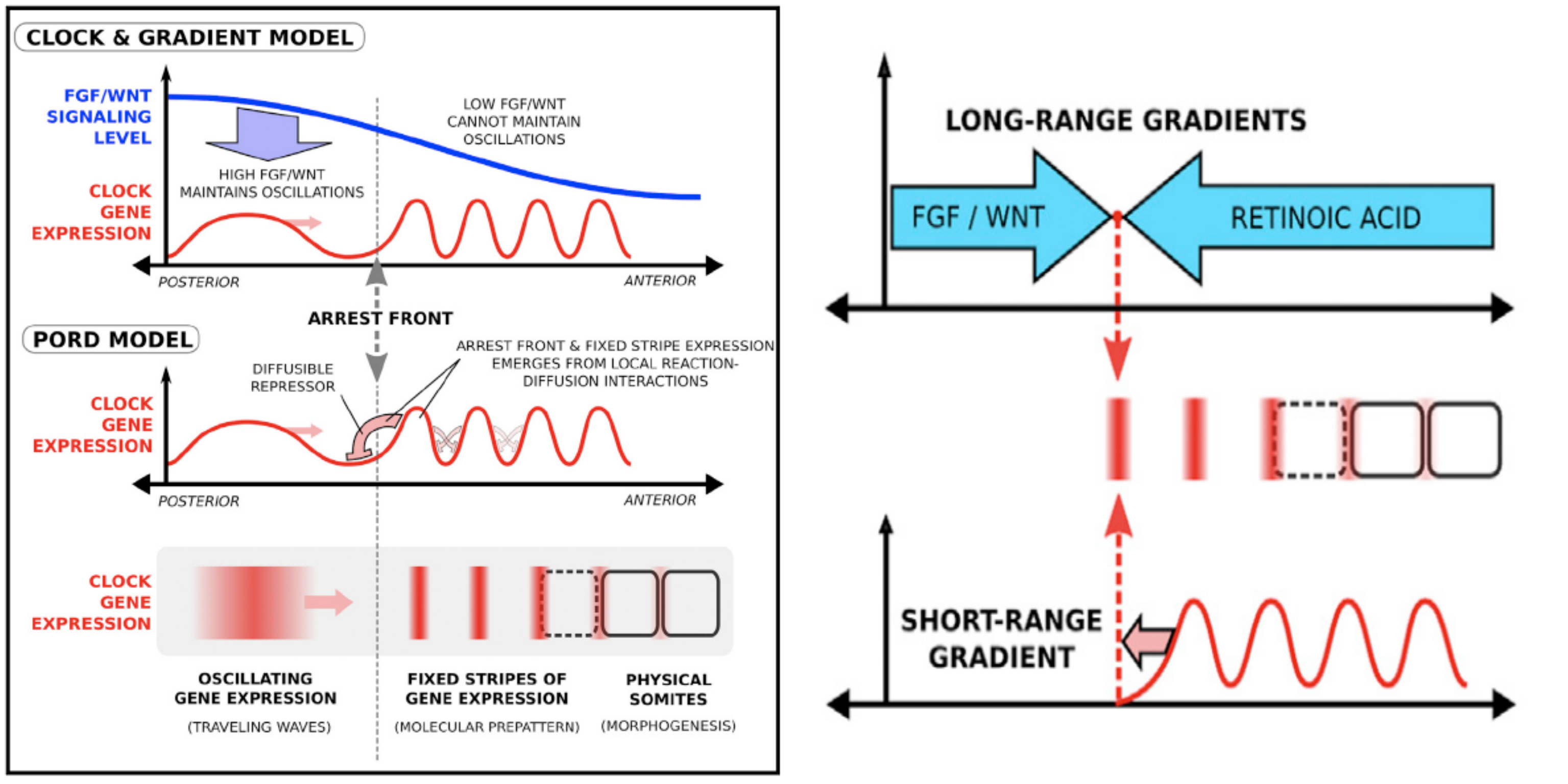}
  \caption{The comparison between C \& W model's and the PORD model's mechanisms\cite{Cotterell}. The left figure shows C \& W model focuses on global gradient control while the PORD model focuses on short-range interactions between cells. The bottom figure shows the oscillation of gene expression in the PSM cells, which are the stripes. Each stripe of gene expression will, in the future, correspond to a subsequent somite boundary. The right figure is the comparison of the sensitivity of the stripe positions. The positional accuracy of the arrest front will be more sensitive to noise if defined by long-range gradients(top) than if defined by the distance from the last-formed expression stripe (bottom).}
  \label{fig: PORD}
\end{figure}

\par The PORD model argues that there is a molecular patterning process that sequentially produces stripes of gene expression along the PSM, resulting in the segmentations in the PSM. Figure \ref{fig: PORD} shows this mechanism and its comparison with the C \& W model. There are two dynamical systems that are involved in this process. First, cells of the PSM exhibit oscillations of gene expression. These oscillations are organized into traveling waves and they are locally well synchronized: neighboring cells are in very similar phases of the cycle \cite{Cotterell}. Second, these oscillations are arrested in an anterior-to-posterior progression, which means the position where the oscillations are frozen travels posteriorly through the PSM, and that position is called the arrest front. Note that the arrest front is similar but not equivalent to the determination front mentioned above, and is addressed in the discussion section in Cotterell's paper \cite{Cotterell}. 
\par Despite being locally self-organizing, the PORD model involves both molecular oscillations in the PSM and a traveling wavefront. Yet, it continues to create stripes even in the absence of a moving FGF gradient. Thus it does not rely on positional information along the PSM. In this reaction-diffusion model, the distance between stripes is defined by the local diffusion of a repressor molecule, which is secreted from the stripes themselves (See the right of Figure \ref{fig: PORD}). However, the fact that the model behaves the same with and without the gradient seems like a potential problem, since the PSM cells have been studied without a gradient and their behaviors seem to be very different \cite{Hubaud}.
\par Overall, the PORD model challenges the existing clock and wavefront models by providing a fundamentally alternative theory based on locally self-organization. It could explain somite size scaling and have higher robustness of somite size regulation. Some of the PORD model's results also stand the test in chick embryos, which shows its validity\cite{Cotterell}. 

\subsection{Mathematical Equation}
\par The exploration of mathematical equations for the PORD model is refreshing. Cotterell and colleagues enumerated all possible topologies that are possible for a gene regulatory network of three genes. Of the 9710 possible networks, 210 produced a multi-stripe pattern for at least one parameter set. Of all the stalactites in the topological tree containing successful topologies, they found two versions of the C \& W model and several versions of the oscillatory PORD model. 
\par
\begin{figure} [H]
  \centering
  \includegraphics[width = \linewidth]{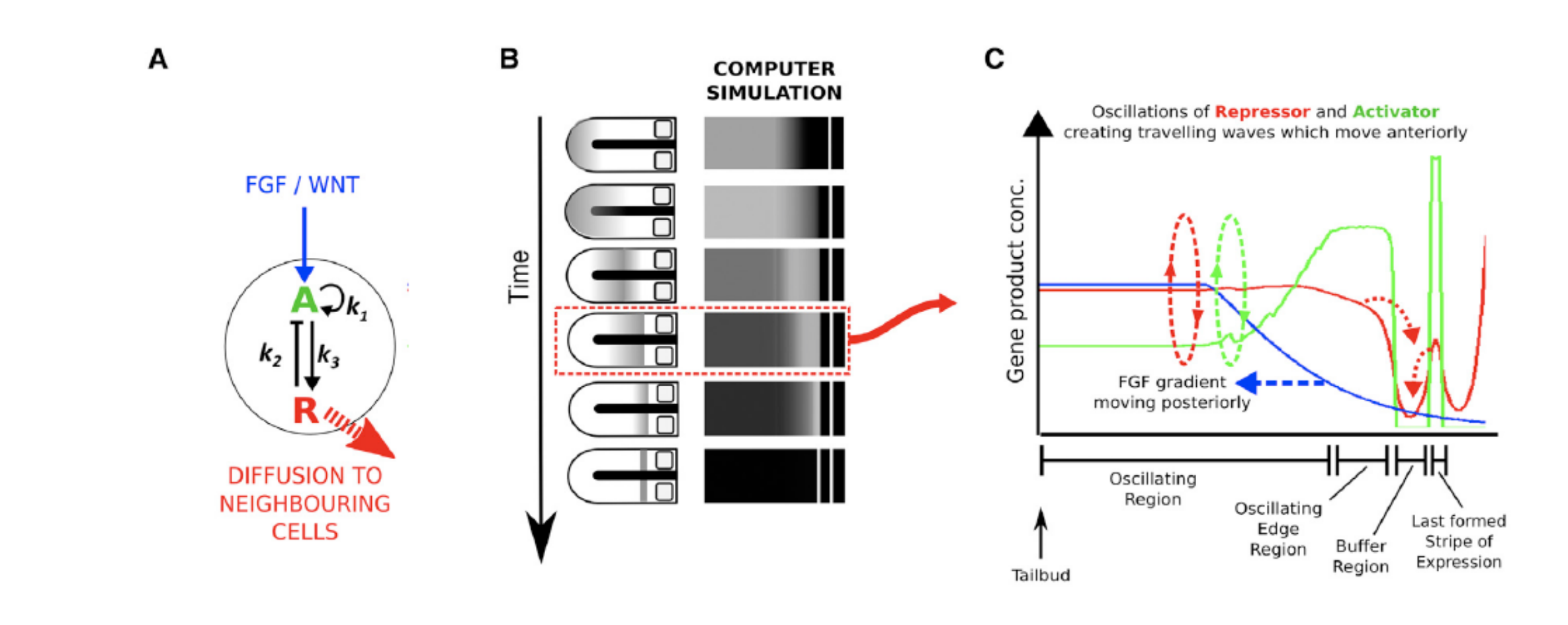}
  \caption{The PORD mechanisms \cite{Cotterell}. \textbf{(A) }The minimum somite-patterning circuit that implements the PORD mechanism. It contains an activator molecule (green) and a diffusible repressor (red). $\kappa_1, \kappa_2,$ and $\kappa_3$ are strengths of interactions between A and R. It shows that the system relies on the FGF to activate but does not need FGF to perpetuate, which fits the PORD character. \textbf{(B)} Gene expressions are initiated at the posterior end and they travel to the anterior end. The white stripe is the formed somite, which is the position where the last gene expression stopped. When the next wave propagation arrives at a certain distance from the last one, it will stop and form the next stripe. \textbf{(C)} A snapshot of gene expression oscillation along the PSM. The oval dashed arrows indicate the oscillation directions. The blue line is the FGF gradient. It's not directly related to stripe formation. The Buffer Region is generated by the diffusion of the repressor from the last formed stripe, which inhibits oscillations. Therefore, cells cannot go beyond the Buffer Region and will exit oscillations in the Oscillating Edge Region to form new stripes. Newly formed stripes act as the next source of repressor to prevent oscillations. They will form new Buffer Regions and push the arrest front posteriorly.}
  \label{fig: PORD2}
\end{figure}

 The simplest design of the oscillatory model is a network that contains only two nodes ((A) of Figure \ref{fig: PORD2}), comprising a cell-autonomous activator (A), which is itself activated by the FGF signal, and a diffusible repressor (R). A and R are defined by the following equations: 
\begin{gather}
\frac{\partial A}{\partial t} = \Phi (\frac{\kappa_1 A - \kappa_2 R +\textit{F} + \beta}{1 + \kappa_1 A + \kappa_2 R + \textit{F }+ \beta}) - \mu A \\
\frac{\partial R}{\partial t} =\frac{\kappa_3 A}{1+\kappa_3 A} - D \nabla^2 R - \mu R
\end{gather}
where $\kappa_1, \kappa_2, \kappa_3$ define the strengths of regulatory interactions between A and R. D is the diffusion constant for R, $\mu$ is a fixed decay constant, and F is the regulatory input of the FGF gradient onto a. $\beta$ is the background regulatory input of A. To prevent negative values of morphogens, we use the function $\Phi(x) = xH(x)$, where $H(x)$ is the Heaviside function.
\par Together they form a reaction-diffusion mechanism where R inhibition is responsible for the spacing of adjacent stripes. Since the PORD model does not rely on global positional information, the model does not spontaneously generate segments everywhere but rather progresses from anterior to posterior, which is similar to real-world biological phenomena. (B) and (C) in Figure \ref{fig: PORD2} shows the wave of gene propagation and its oscillatory mechanism.

\subsection{Analysis}
\par
\begin{figure} [H]
  \centering
  \includegraphics[width = \linewidth]{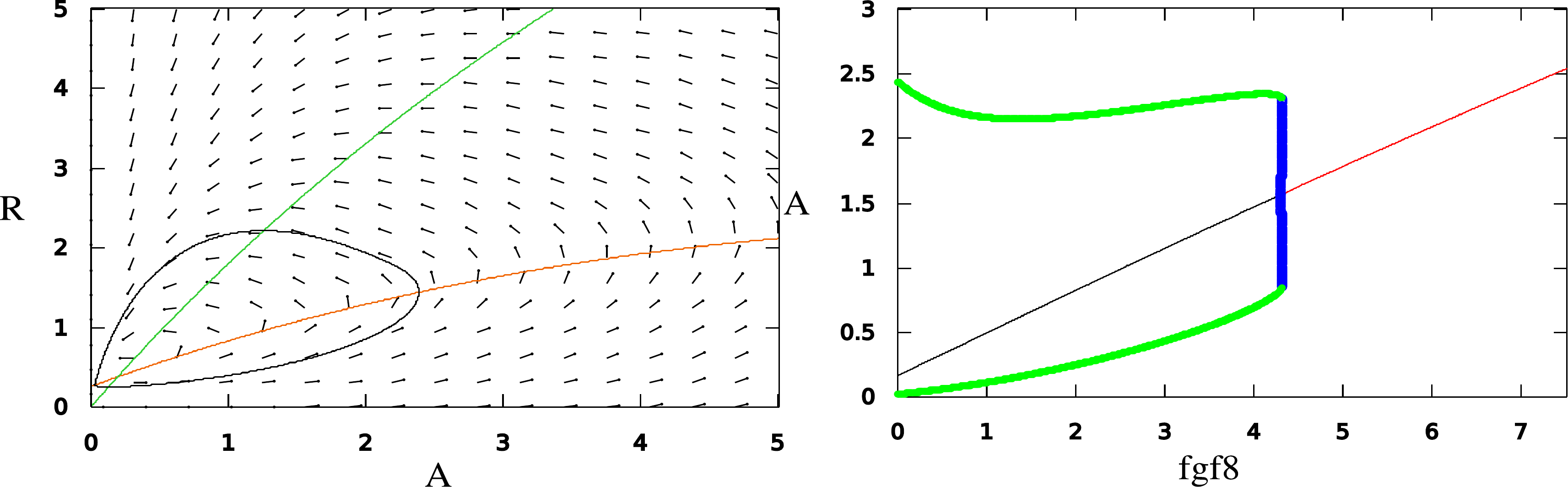}
  \caption{XPP analysis of the PORD model. The left figure is the phase portrait for the non-diffusing case of Equations (8) and (9), which means to ignore the diffusing term D. The green and red lines are the nullclines for the activator A and inhibitor R respectively. The right figure is the bifurcation analysis for the same case. It shows the bifurcation of activator A with respect to the FGF gradient.}
  \label{fig: XPP2}
\end{figure}

 The PORD model proves to be a typical oscillatory model. Its wave propagation theory as well as the mathematical equation both exhibit its oscillatory nature. We used XPP and recapitulated figure E of Figure 2 in the Cotterell paper \cite{Cotterell}, see the left of Figure \ref{fig: XPP2}. 

Ignored its diffusion state and made stationary, the system reveals that oscillations are the natural dynamic state for most cells in the PSM. The bifurcation analysis (the right of Figure \ref{fig: XPP2}) also reveals its oscillatory nature. The activator has a Hopf bifurcation which starts when the fgf8 drops to a certain level. Meaning, when the fgf8 is high, the activator will be stimulated. When fgf8 decreases to a certain level, the activator will interact with the inhibitor and will start to oscillate. The part in between the two green boundaries is where oscillation exists. And the oscillation will stop once the cells reach the arrest front, which in this case, is when the fgf8 decreases to 0. 
\par However, the PORD model has received some criticisms. For one, although in the paper, it's been claimed multiple times that this model does not require the moving FGF gradient, it nevertheless acts to couple the rate of embryo growth with the integral levels of FGF8 signaling in the PSM \cite{Cotterell}. Meaning, it can't ignore the fact that FGF8 plays an important role in controlling the size of somites. Higher levels of FGF signaling will result in smaller somites. Also, from the bifurcation analysis shown above, we can see that the character of the PORD system is somewhat similar to C \& W system in that FGF8 gradient information could control the cell activities in both cases. The position where the activator starts to oscillate can be seen as the determination front in C \& W model, and the new terminology "arrest front" in the PORD system also locates close to where the FGF8 gradient drops to a very low level. Simply put, although the PORD system has created some new terms such as the "buffer region" and "the arrest front", its activities, similar to C \& W model, could still be explained by the FGF8 gradient control. Another problem is that, not only in PORD model but also in some other oscillatory-based models, many of them control the specific position of spatial stripes "manually", by defining thresholds or piecewise functions \cite{Fhn}. Although this may help create beautiful results, it's in contrast to the principle of the self-organization of biological systems. 
\par Nonetheless, the PORD model, as one of the most famous oscillatory models for somitogenesis, does present a very different perspective. It reveals the possibility that cells themselves carry an oscillatory nature in the absence of diffusion. It also produces several nice movies to show the oscillation process clearly. We tried to reproduce that in Matlab but didn't succeed. 

\section{Excitable Model}
\subsection{Summary of the one-dimensional RD Model}
Both the C \& W model and the PORD model have an important feature that has not been mentioned in previous sections. That is, both models set spatial continuity as a key requirement. Spatial continuity, in this case, means that both models ignore the size of cells and see the PSM as a whole unit, or, a spatial continuum. However, although spatial continuity is acceptable in most chemical reaction systems, Nagahara and colleagues argue that it's not always the case in biological systems. It's simply because cells in a multicellular organism have a finite size\cite{Fhn}. In the initial stages of an organism's developmental process, when important biological structures first emerge, the number of cells is usually small and the size of a single cell can not be simply ignored since the size of the field where the phenomena occur is comparable to that of a cell. As we know spatial continuity is not always met, we may have to consider the spatial variations between cells. Yet, Nagahara and colleagues propose that instead of thinking in that way, it's better to just treat cells as "interacting discrete nodes" in a network. Since the diffusion inside a cell is much faster than that in a membrane, treating cells as individual nodes that compose a huge network is suitable. 
\par
\begin{figure} [H]
  \centering
  \includegraphics[scale = 0.3]{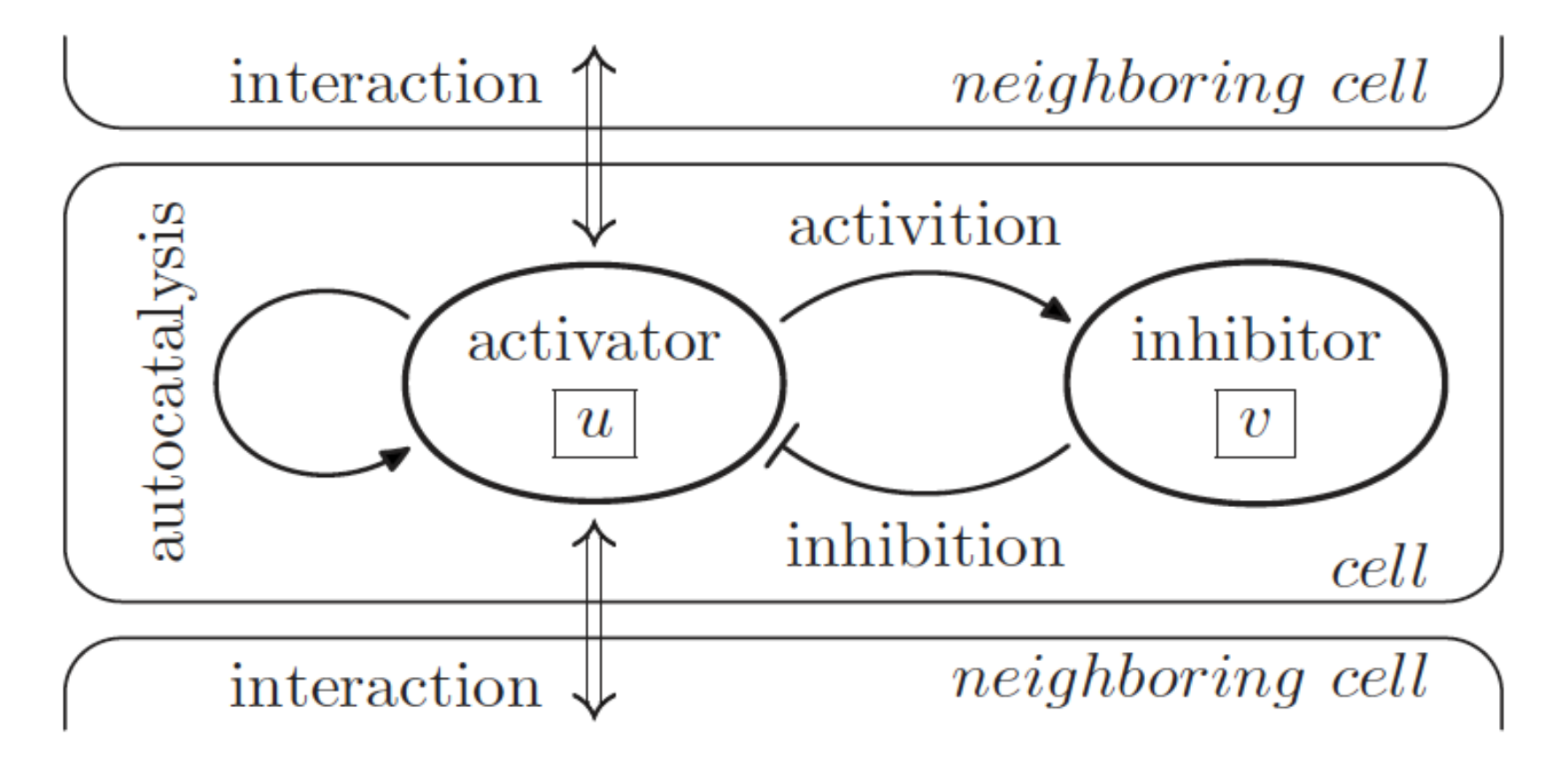}
  \caption{An illustration of the one-dimensional reaction-diffusion system proposed by Nagahara and colleagues \cite{Fhn}. The major difference between this one and the PORD model is that the diffusion of the inhibitor v is excluded. There's only one interaction between neighboring cells, which is the activator u.}
  \label{fig: Fhn}
\end{figure}

 Nagahara also criticized the complicity and difficulty in other models, since models based on a continuum will have difficulty producing a narrow boundary between different distinctive behaviors, as sharp as two or three cells. The fact that most models assume the existence of two or more different interactions, such as activator, inhibitor, etc., among the neighboring cells, is also complicated. Therefore, they created a simple, one-dimensional reaction-diffusion model that focuses on three things: (1) No diffusion among inhibitors, (2) cells are discrete instead of continuous, and (3) spatial inhomogeneity, where (1), (2) are new ideas, while (3) is an old one. See Figure \ref{fig: Fhn}.

\subsection{FhN-type system and excitability}
\par Below is a hypothetical model that describes the above features of gene expression in somitogenesis: 
\begin{gather}
\frac{\partial u}{\partial t} = f(u,v) + D \mathcal{L} u \\
\frac{\partial v}{\partial t} =g(u,v)
\end{gather}
where u and v are concentrations of the activator and inhibitor respectively. $D \mathcal{L}u$ is the diffusion term for the activator while f and g are the reaction terms of the activator and the inhibitor. f and g are given as follows:
\begin{gather}
f(u,v) = \frac{1}{\tau_1}(\frac{1}{\gamma}u(u-a)(1-u)-v+\beta) \\
g(u,v) = \frac{1}{\tau_2}(u-v)
\end{gather}
where $\tau_1$ and $\tau_2$ represent the time scales of local reaction kinetics of u and v respectively. $\gamma$ represents the spatial gradient and the temporal change in the concentration of a certain substance \cite{Fhn}, which is dependent on space x and time t. In reality,$\gamma$'s biological counterpart, in the case of somitogenesis, could be the FGF8 gradient in the PSM. $\gamma$ thus plays an important role in this model.
\par We call this model the "FhN-type" model. The reason for that is because the reaction terms (12) and (13) resemble the "Fitzhugh-Nagumo" model. It's named after Richard Fitzhugh who suggested the model in 1961 and J. Nagumo who created the equivalent circuit the following year. The Fitzhugh-Nagumo model is a generic model for \textit{excitable systems}. Because of its simple two-variable form and generality, it has been used widely. The Fitzhugh-Nagumo prototype model has the following form:
\begin{gather}
\dot{v} = v-\frac{v^3}{3}-w+I_{ext} \\
\tau \dot{w} = v+a-bw
\end{gather}
\par
\begin{figure} [H]
  \centering
  \includegraphics[scale = 0.3]{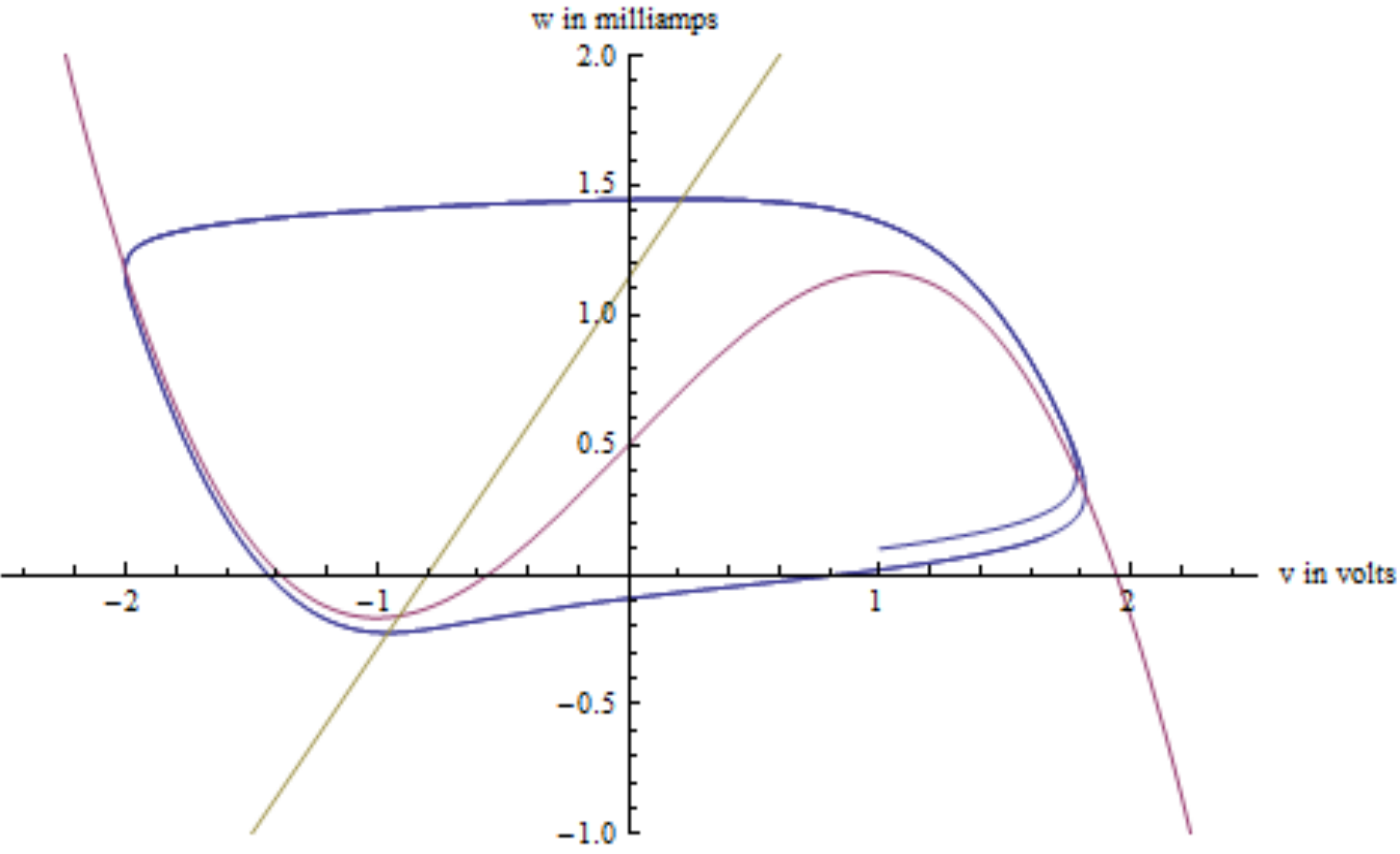}
  \caption{The phase plane for the prototype Fitzhugh-Nagumo model. $I_{ext}$ = 0.5, a = 0.8, b = 0.7. The blue line is the trajectory of the FhN model in the phase space. The pink and yellow lines are two nullclines where the pink line is a cubic nullcline and the yellow line is the linear nullcline.}
  \label{fig: excitable}
\end{figure}

 The reason why we say this model represents \textit{excitable systems}, is because when $I_{ext}$, the external stimulus, exceeds a certain threshold level, the system will exhibit a characteristic excursion in the phase plane, before the variables $v$ and $w$ relax back to the rest values. We say the system is excited and will be refractory to excitability for a period of time. When $I_{ext}$ doesn't exceed that threshold, there will not be an excursion, and we say the system is not yet excited and remains quiescent or excitable. Except for the excursion, the phase plane also contains two nullclines. One is linear and the other is a cubic nullcline, or a sigmoid. The excitability of the system can be discovered by looking at the spatial relationship between the two nullclines: the closer the linear nullcline to the peak of the sigmoid, the more excitable the system is. See Figure \ref{fig: excitable} for the phase plane.

\subsection{Analysis}
\par
\begin{figure} [H]
  \centering
  \includegraphics[width = \linewidth]{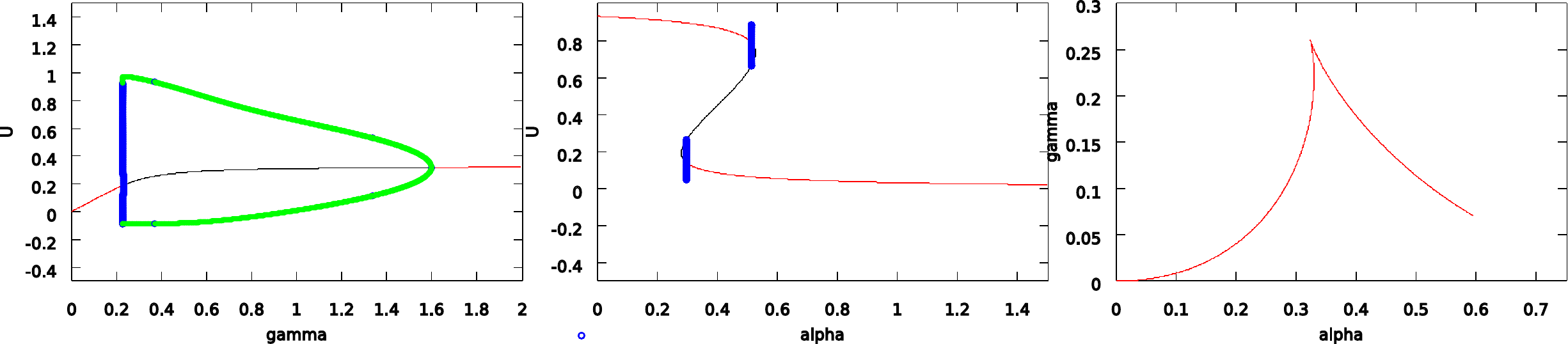}
  \caption{The bifurcation diagrams for the one-dimensional RD model. The left figure is a hopf bifurcation of $\gamma$. The middle figure is a saddle-node bifurcation of $\alpha$. The right figure is a co-dimension bifurcation diagram, which is a cusp bifurcation of $\gamma$ and $\alpha$. Cusp bifurcation is defined as two branches of saddle-node bifurcation curve that meet tangentially, forming a semi-cubic parabola.}
  \label{fig: Fhnbif}
\end{figure}

 In this model, $\gamma$ is defined as a linear gradient function: 
\begin{gather}
\gamma(x) = 0.21-0.20x, \quad \quad \quad x \in [0,1]
\end{gather}
With the given information, we used XPP and created the bifurcation diagrams to find the model's qualitative features. See Figure \ref{fig: Fhnbif}. We can see that $\gamma$ has a hopf bifurcation. When $x \in [0,1]$, $\gamma$ is invariant and the system is stable, yet when x increases beyond the boundary, the system tends to become oscillatory. $\alpha$'s bifurcation is a saddle-node bifurcation. Together with the stable region of $\gamma$'s bifurcation, the system tends to exhibit a bistable region when $\gamma$ is small. Two stable regions will coexist and then the system will tend to be oscillatory as $\gamma$ increases. The right figure is a cusp bifurcation. Note that the cusp bifurcation has a normal form which resembles equation (12). 
\par
\begin{figure} [H]
  \centering
  \includegraphics[scale = 0.25]{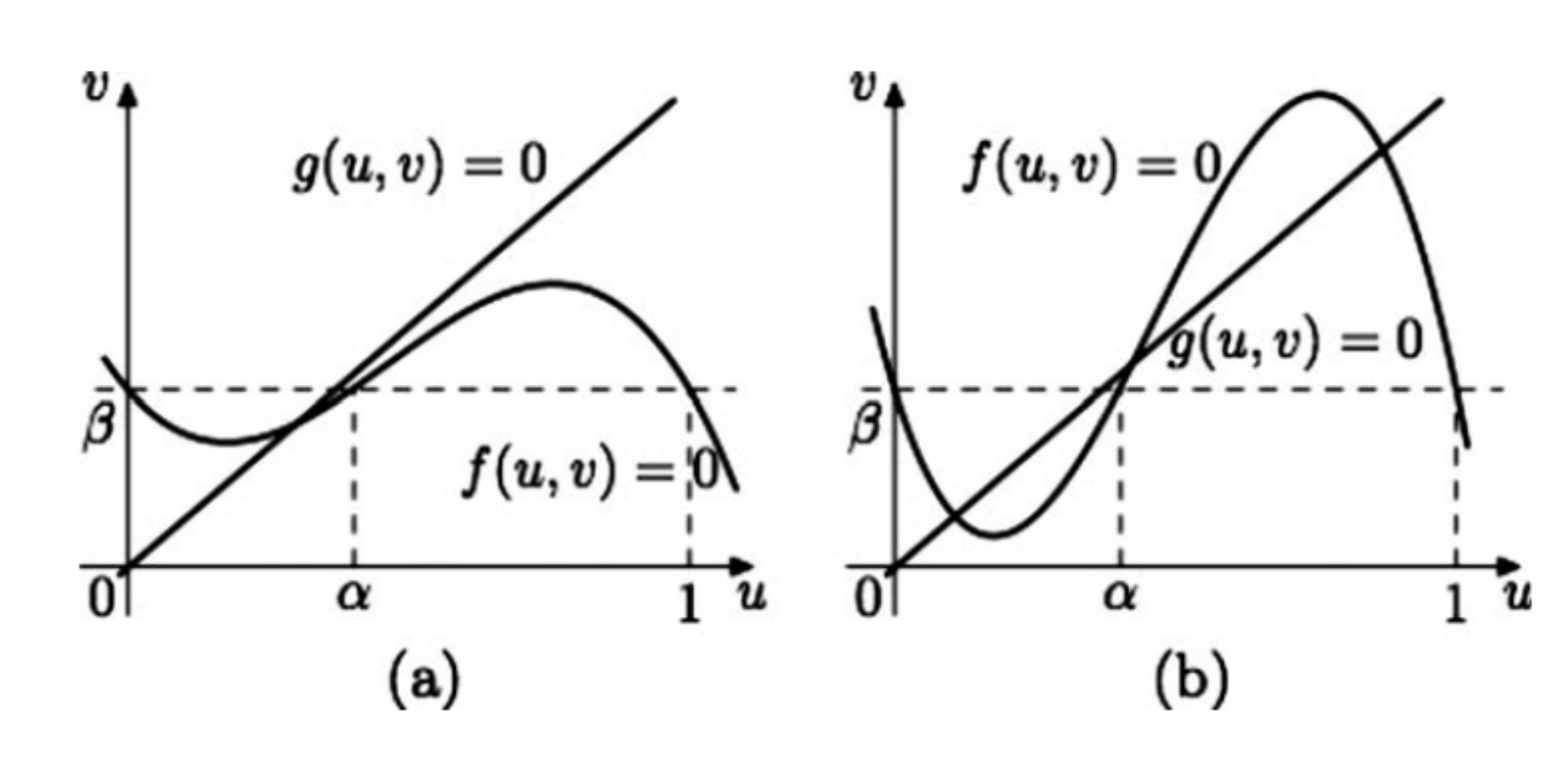}
  \caption{Phase planes for equation (12) and (13)\cite{Fhn}. They correspond to (a) oscillatory state in big $\gamma$ case and (b) bistable state in small $\gamma$ case.}
  \label{fig: phase}
\end{figure}

\begin{figure} [H]
  \centering
  \includegraphics[width = \linewidth]{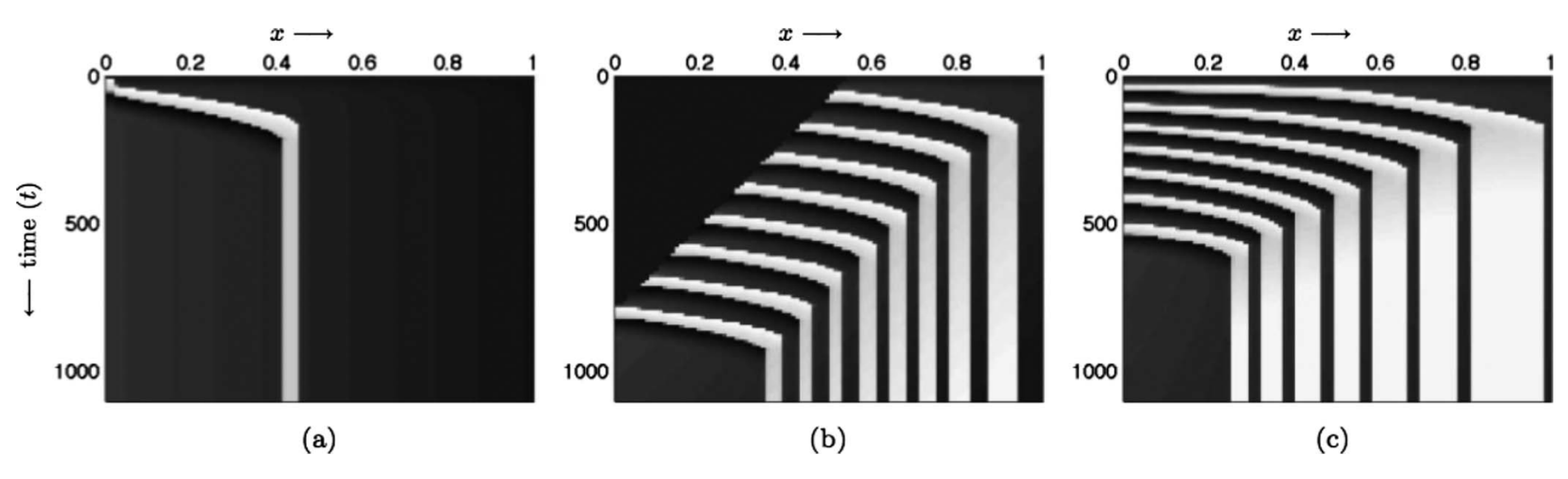}
  \caption{The numerical simulations for the one-dimensional RD model, illustrate the manner of wave propagation, depending on the geometrical distribution of $\gamma$ \cite{Fhn}. }
  \label{fig: result2}
\end{figure}

 The Nagahara paper utilizes the idea of bistability. By changing the parameter $\gamma$, which adjusts the amplitude of cubic function, we can vary the local kinetics from oscillation to bistability \cite{Fhn}. Figure \ref{fig: phase} shows the changes in the nullclines $f(u,v) = 0$ and $g(u,v) = 0$ by manipulating $\gamma$. As we can see, the system gains more excitability when $\gamma$ is decreased, since the decrease in $\gamma$ will make the sigmoid, or the cubic nullcline, elongated vertically. As the linear nullcline will not change, this elongation will close the distance between the linear nullcline and the sigmoid, thus making the system more easily excited. Furthermore, using plotting techniques in Mathematica, we found that manipulating $\alpha$ and $\beta$ in equation (12) will also change the excitability, in which the increase in $\alpha$ makes the system less excitable, while the increase in $\beta$ makes the system more excitable. 
A spatiotemporal diagram is shown in Figure \ref{fig: result2}, where (a) $\gamma$ is given as invariant. A single pulse triggered from the left boundary propagated to the right and generated a stationary band at a specific position. (b) If we take into account the posterior growth of the PSM, $\gamma$ will be a function of both space and time. The pulses will lead to a static, periodic structure. (c) $\gamma$ sets to decrease as the wave passes, while we ignore the growth of the PSM. The pulses can also create a static, periodic structure, but the bandwidth will be much thicker. 
\par However, the Nagahara model is not a typical RD model. The structures in Figure \ref{fig: result2} will not form unless the model is discrete, which is in opposition to continuous. Spatial discreteness is an important feature, and also the base for this model. Normally in the continuous case, the wave propagation triggered from the left will not stop and generate a stationary band, rather, it will propagate across the field without stopping\cite{Fhn}. The paper proposes that, when D is small enough, the propagation of the wave will be blocked and there exist stable steady solutions. This phenomenon is called "wave propagation failure". We used Matlab and simulated what will happen to the system as D varies, and we obtained the result shown in Figure \ref{fig: failure}, using the code provided in the appendix. 
\begin{figure} [H]
  \centering
  \includegraphics[width = \linewidth]{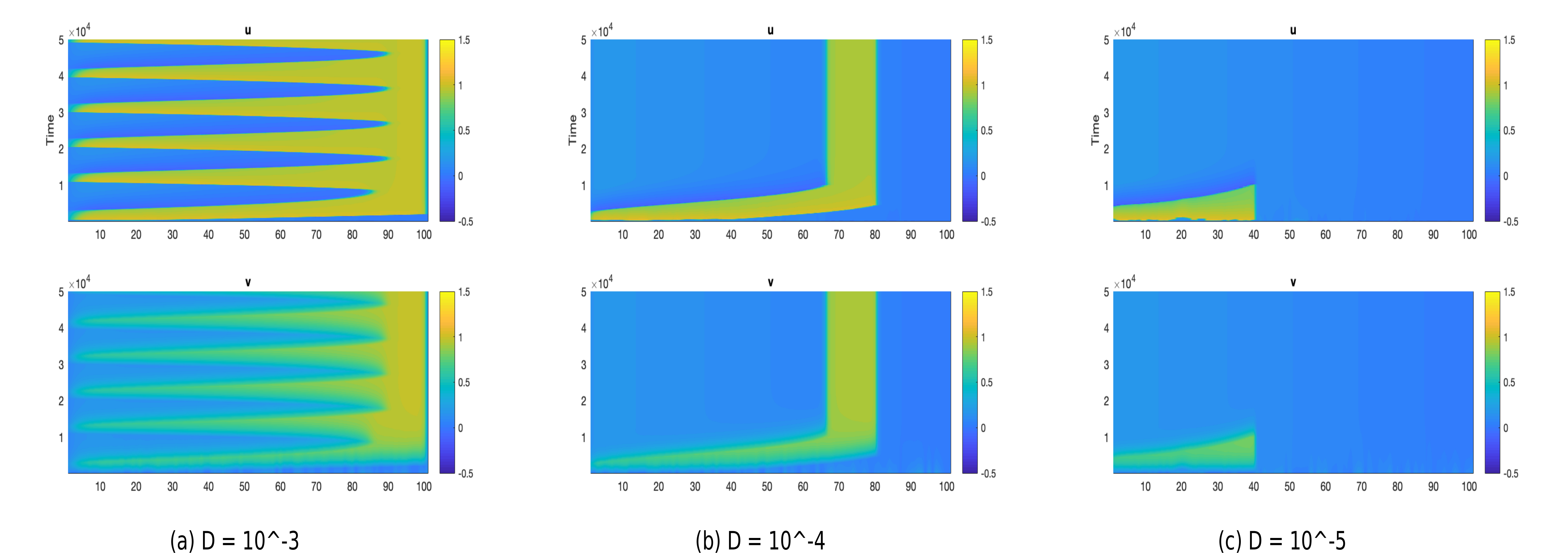}
  \caption{The spatiotemporal diagram for the system as D varies. D decreases from (a) to (c). In (a), $D = 10^{-3}$. The system is not yet discrete and the wave propagates through the field non-stop. In (b), $D = 10^{-4}$. The system becomes discrete, and the wave propagation leaves a steady state solution, just like what's shown in Fig \ref{fig: result2} (a). In (c), $D = 10^{-5}$. The system is discrete but the diffusion is so weak that the signal fails to reach the bistable region.}
  \label{fig: failure}
\end{figure}
\par Overall, this one-dimensional RD model, which is the FhN-type model, is the most immature model among the three. Unlike the previous two that have already been going through in vitro experiments, this model is highly theoretical and leans more toward physics than biology. However, the simple, fresh idea of excitability opens a new perspective to see the whole process. Excitability in this model, as we have found in Mathematica, can vary with respect to $\gamma, \alpha$, and $\beta$. Depending on three variables seems a bit capricious but the simple philosophy behind excitability - how easy it is to cross the threshold - makes it simple to accurately manipulate the excitability of a system or how to make it "excited". We look forward to exploring more about this feature, and to implementing it in an effort to effectively understand somitogenesis.

\label{sec: excitable}
\label{sec: other}

\section{Conclusion}
\label{sec: conclusion}
\par As they represent several of the mainstream ideas in the field, it is not surprising that all three models provide abundant insights into the hidden mechanisms behind somitogenesis. While their core ideas are different in one way or another, and each of them has its own flaws, each model's existence nevertheless greatly improves scientists' understanding toward this field and motivates new experiments. 
\par The Clock and Wavefront model proposed a prescribed determination front, which is determined by the level of FGF8 gradient, that controls the positioning of somites. It segments the PSM cell into different regions and explains the somitogenesis process systematically. The result is easy to recapitulate and it shows desired characters that fit the theory. Its biggest drawback is its hard-coded outcomes, for example, the model can not explain experiments in which mutant embryo determination fronts change since its position is prescribed in the model. 
\par The PORD model proposes a local reaction-diffusion oscillatory mechanism that can generate stripes of somites without global gradient information. The system is simple and their video result is very impressive. However, in contrast to its claims, the PORD model still relies on the FGF8 gradient to control the size of somites, and many of its mechanisms can still be explained by \textit{fgf8} which seems as if the theory is like another "perspective" to see the \textit{fgf8}'s effects. Many of the oscillatory models control the specific position of the stripes manually using thresholds or piecewise functions, which violates its self-organizing nature. 
\par The excitable model is a one-dimensional reaction-diffusion model that resembles the "Fitzhugh-Nagumo" model, which is an excitable, generic model widely applied in fields of physics. Discreteness plays an important role in the system: the model decreases the diffusion level until cells are no longer continuous, then blocks the wave propagation generated by excitability, creating fixed stripes of somites, which is a refreshing idea. This model is theoretical and has not gone through in vitro experiments, and its result is not robust as it needs a very fine-tuned diffusion level, which is hard to achieve in reality, but the idea of an excitable system has latent potential and much to be exploited. 
\par The C\&W and PORD models yield great results that match our desire. However, their strict conditions often do not meet in real-life biology, and it seems that this problem can't simply be solved under the existing frames of mathematical equations.  In my perspective, excitability is what needs to be studied the most rigorously in order to understand somitogenesis. Dr. Hubaud and his colleagues' research \cite{Hubaud} which proposed excitability as a general framework for oscillations in the PSM cell is a great start. While referring to other models for inspiration, we should boldly explore the direction of excitable models, instead of sticking to past experiences and techniques on other models that seem more successful at the moment, since forsaking the mindset that holds for the moment is a must to create a better one.

\bibliographystyle{plain}
\bibliography{refs}

\begin{appendices}
\section{Code}
\subsection{Clock and Wavefront}
\begin{lstlisting}
function newclockwaves
m = 0; % a slab system (in terms of slab,cylindrical,spherical)
x = -5:0.01:15; % setting up AP axis
t = 0:0.05:25; % setting up time axis

sol = pdepe(m,@pdex4pde,@pdex4ic,@pdex4bc,x,t);
u = sol(:,:,1);
v = sol(:,:,2);
w = sol(:,:,3);

figure
imagesc(x,flipud(t),u);
colormap(gray)
colorbar
set(gca,'YDir','normal')
title('(a) Somitic Factor')
xlabel('AP axis (x)')
ylabel('Time (t)')

figure
imagesc(x,flipud(t),v);

colormap(gray)
colorbar
set(gca,'YDir','normal')
title('(b) Signaling Molecule')
xlabel('AP axis (x)')
ylabel('Time (t)')

figure
imagesc(x,flipud(t),w);
colormap(gray)
colorbar
set(gca,'YDir','normal')
title('(c) FGF8')
xlabel('AP axis (x)')
ylabel('Time (t)')
% --------------------------------------------------------------
function [c,f,s] = pdex4pde(x,t,u,DuDx)
x1 = 1;
x2 = 0;
xn = 0;
k = 10;
mu = 0.0001;
epsilon = 0.001;
gamma = 0.001;
Dv = 50;
Dw = 20;
n = 1;
o = 0;
cn = 0.5;
xb = 7.5;
epsilon1 = 0;
%n1=(-cn-sqrt(cn^2+4*n*Dw))/(2*Dw);
%n2=(-cn+sqrt(cn^2+4*n*Dw))/(2*Dw);

c = [1; 1; 1]; 
f = [0.00001; Dv; Dw] .* DuDx; % these are the diffusion coefficients, take D1 \approx 0 for this model
%w1 = n1/(n*(n1-n2))*exp(n2*(-xn+1));
Xu = heaviside(x1-x+cn.*t);
Xv = heaviside(x2-x+cn.*t);
Xw = heaviside(x-xn-cn.*t);
Xb = heaviside(epsilon1-xb+x)*heaviside(epsilon1+xb-x);

F1 = ((u(1)+mu*u(2)).^2)./(gamma+u(1).^2).*Xu-u(1); 
F2 = k.*(Xv./(epsilon+u(1))- u(2));
F3 = Xw+o*Xb-n.*u(3);

s = [F1; F2; F3]; 
%% Set the initial conditions
function u0 = pdex4ic(x)
Dv = 50;
Dw = 20;
%gamma = 0.001;
k = 10;
xn = 0;
n = 1;
cn = 0.5;
epsilon = 0.001;
epsilon1 = 0;

lam = sqrt(k/Dv);
A = 1/(1+epsilon-epsilon1);
B = A*sign(x)/(2*cosh(lam*10));
%n0 = 1/2*(1+sqrt(1-4*gamma/k));
n1=(-cn-sqrt(cn^2+4*n*Dw))/(2*Dw);
n2=(-cn+sqrt(cn^2+4*n*Dw))/(2*Dw);
w0=heaviside(xn-x)*(n1/(n*(n1-n2))*exp(n2*(x-xn)))+
heaviside(x-xn)*(n2/(n*(n1-n2))*exp(n1*(x-xn))+1/n);
u0 = [heaviside(-x); A*heaviside(-x)+B*cosh(lam*(10-abs(x))); w0]; 
%% Set the boundary conditions
function [pl,ql,pr,qr] = pdex4bc(xl,ul,xr,ur,t)
epsilon = 0.001;
gamma = 0.001;
k = 10;
Dv = 50;
Dw = 20;
n = 1;

%pl = [ul(1).^2./(gamma+ul(1).^2)-ul(1); k.*(1/(epsilon+ul(1))-ul(2)); 0];
%ql = [0.00001; Dv; Dw]; % can't be set to [0,D], need to use a small value to approximate 0
%pr = [0; 0; 1-n.*ur(3)];
%qr = [0.00001; Dv; Dw]; 
pl = [0;0;0];
ql = [1;1;1];
pr = [0;0;0];
qr = [1;1;1];
\end{lstlisting}

\subsection{Nagahara Discrete}
\begin{lstlisting}
%% 1D Simulation of spatial savanna model with diffusion in space
tic
%% Set options for plots/movies
close all
fprintf('\n');
DO_MOVIE = 1; % i.e. write movie to avi file for playback, else just plots
ONE_D_PLOT = 0; % plot dynamics of a single point over time
END_STATE_PLOT = 0; % 2D plot of the final state of the system
THREE_D_PLOT = 0; % plot two spatial dimensions, time and colour
SPACE_TIME_PLOT=1;
%% Numerical method parameters
L = 1; % working on [0,L]
N = 100; % N+1 grid points
delta = L/N;  % spatial discretization parameter -> 0.01 suggested
h = 0.01; % time discretisation parameter
n = 50000; % number of time steps
tau = (n-1)*h; % simulations time domain is [0,tau]
%% Function definitions
gamma_fun = @(x, a, b) a + b.*x;
f_fun = @(u,v,tau1,gamma,alpha,beta) (u.*(u-alpha).*(1-u)./gamma-v+beta)/tau1;
g_fun = @(u,v,tau2) (u-v)/tau2;
%% Model parameters (values from Nagahara et al. (2009), Phys. Rev. E)
tau1 = 0.588;
tau2 = 32.1;
alpha = 0.4;
beta = 0.33;
a = 0.21;
b = -0.2; % b = -0.2 in the paper
D1 = 1*10^(-5); % D too big or small won't work
D2 = 0; % this is set to zero in Nagahara's paper
%% Set up the initial distributions of the cover types on the grid
% each row is one time step of the simulation
% the initial condition is a pulse of u near x = 0
u(1,:) = 0.2*rand(1,N+1);
v(1,:) = 0.2*rand(1,N+1);
% Compute the birth and mortality matrices as a function of rainfall
X = 0:delta:L;
gamma_grad = gamma_fun(X,a,b); % compute the rainfall gradient along the x-axis
%% The numerical scheme
for i = 2:n
    % compute convolutions for this time step
    progressbar(i,n);
    u(i,:) = u(i-1,:) + h*( f_fun(u(i-1,:),v(i-1,:),tau1,gamma_grad,alpha,beta) + ...
        D1*([0, u(i-1,1:(end-1))] - 2*u(i-1,:) + [u(i-1,2:end),0])./(delta*delta) );
    % NB the zeros reflect that we are using an "open boundary"
    v(i,:) = v(i-1,:) + h*( g_fun(u(i-1,:),v(i-1,:),tau2) + ...
        D2*([0, v(i-1,1:(end-1))] - 2*v(i-1,:) + [v(i-1,2:end),0])./(delta*delta) );
end
% The following output is useful when trying to discern whether or not 
% a solution is stationary in time 
fprintf('\n');
fprintf('The maximum changes on the grid for each variable at the last time step were:\n');
fprintf(['u: ',num2str(max(abs(u(n,:)-u(n-1,:)))),'\n']);
fprintf(['v: ',num2str(max(abs(v(n,:)-v(n-1,:)))),'\n']);
toc
%% Visualise the solution...
% either in a movie...
if DO_MOVIE
    subplot_vis = struct('cdata',[],'colormap',[]);
    w = VideoWriter('reaction_diffusion_model.avi');
    w.FrameRate = 4;
    open(w)
    f = figure;
    for j = 1:300:n
        f.Name = ['Simulation time: t = ', num2str((j-1)*h)];
        ax1 = subplot(2,1,1);
        plot(X,u(j,:));
        title('u');
        
        ax2 = subplot(2,1,2);
        plot(X,v(j,:));
        title('v');
        
        xlim([ax1 ax2],[0 L])
        ylim([ax1 ax2],[-0.5 1.5])
        writeVideo(w,getframe(gcf))
        subplot_vis(j) = getframe(gcf);
    end
    close(w)
end
%% Space-time plot of dynamics
if SPACE_TIME_PLOT
    fST = figure;
    fST.Name = 'Evolution over time';
    
    subplot(2,1,1)
    h1 = pcolor(u);
    shading interp
    title('u');
    colorbar
    set(h1, 'EdgeColor', 'none');
    ylabel('Time');
    caxis([-0.5,1.5])
    
    subplot(2,1,2)
    h2 = pcolor(v);
    shading interp
    title('v');
    colorbar
    set(h2, 'EdgeColor', 'none');
    caxis([-0.5,1.5])
end
\end{lstlisting}
\end{appendices}
\end{document}